\ifx\shlhetal\undefinedcontrolsequence\let\shlhetal\relax\fi

\documentclass[11pt]{amsart}

\usepackage{amsmath}
\usepackage{amssymb}

\newtheorem{theorem}{Theorem}[section]
\newtheorem{claim}[theorem]{Claim}

\newtheorem{corollary}[theorem]{Corollary}

\theoremstyle{definition}
\newtheorem{definition}[theorem]{Definition}

\newtheorem{problem}{Problem}[section]

\theoremstyle{remark}

\newcount\skewfactor
\def\mathunderaccent#1#2 {\let\theaccent#1\skewfactor#2
\mathpalette\putaccentunder}
\def\putaccentunder#1#2{\oalign{$#1#2$\crcr\hidewidth
\vbox to.2ex{\hbox{$#1\skew\skewfactor\theaccent{}$}\vss}\hidewidth}}

\def\smallbox#1{\leavevmode\thinspace\hbox{\vrule\vtop{\vbox
   {\hrule\kern1pt\hbox{\vphantom{\tt/}\thinspace{\tt#1}\thinspace}}
   \kern1pt\hrule}\vrule}\thinspace}

\newcommand{\cf}{{\rm cf}}

\def\qedref#1{$\qed_{\reforiginal{#1}}$}

\title[Open and solved problems]{Open and solved problems concerning polarized partition relations}
\author{Shimon Garti}
\address{Institute of Mathematics
 The Hebrew University of Jerusalem
 Jerusalem 91904, Israel}
\email{shimon.garty@mail.huji.ac.il}

\author{Saharon Shelah}
\address{Institute of Mathematics
 The Hebrew University of Jerusalem
 Jerusalem 91904, Israel
 and  Department of Mathematics
 Rutgers University
 New Brunswick, NJ 08854, USA}
\email{shelah@math.huji.ac.il}
\urladdr{http://www.math.rutgers.edu/\char`\~shelah}

\thanks{Research supported by the Israel Science Foundation, grant no. 1053/11, publication number 1012 of the second author}

\begin{document}
\let\labeloriginal\label
\let\reforiginal\ref

\begin{abstract}
We list some open problems, concerning the polarized partition relation. We solve a couple of them, by showing that for every limit non-inaccessible ordinal $\alpha$ there exists a forcing notion $\mathbb{P}$ such that the strong polarized relation $\binom{\aleph_\alpha^+}{\aleph_\alpha} \rightarrow \binom{\aleph_\alpha^+}{\aleph_\alpha}^{1,1}_2$ holds in ${\rm\bf V}^{\mathbb{P}}$. \newline
Nous passons en revue certains probl\`emes non r\'esolus conceront la relation de partition polaris\'ee. Nous en r\'esolvons deux en montrant que pour chaque ordinal limite non inaccessible $\alpha$, il existe un forcing $\mathbb{P}$ tel que ${\rm\bf V}^{\mathbb{P}}$ satisfait la relation polaris\'ee forte $\binom{\aleph_\alpha^+}{\aleph_\alpha} \rightarrow \binom{\aleph_\alpha^+}{\aleph_\alpha}^{1,1}_2$.
\end{abstract}

\maketitle

\newpage

\section{Introduction}

The central notion in this paper is the following:

\begin{definition}
\label{sstrongpolarized} The strong polarized relation. \newline
The strong polarized relation $\binom{\lambda}{\kappa} \rightarrow \binom{\lambda}{\kappa}^{1,1}_2$ means that for every function $c : \lambda \times \kappa \rightarrow 2$ there are $A \subseteq \lambda$ and $B \subseteq \kappa$ such that $|A|=\lambda, |B|=\kappa$ and $c \upharpoonright (A \times B)$ is constant.
\end{definition}

The consistency of the strong polarized relation $\binom{\lambda^+}{\lambda}\rightarrow \binom{\lambda^+}{\lambda}^{1,1}_2$ has been established many years ago by Laver in \cite{MR0371652}, for the specific case of $\lambda=\aleph_0$. Until recent years, it was open whether any uncountable cardinal may satisfy this relation. Moreover, even weaker relations of this kind were unknown to exist on uncountable cardinlas. See, for example, \cite{MR2367118} which asks whether there exists any uncountable cardinal $\kappa$ such that $\binom{\kappa^+}{\kappa} \rightarrow \binom{\kappa^+\ \kappa}{\kappa\ \kappa}^{1,1}_2$.

By previous papers we know how to force the strong relation $\binom{\lambda^+}{\lambda}\rightarrow \binom{\lambda^+}{\lambda}^{1,1}_2$ for every supercompact cardinal $\lambda$ (see \cite{MR3201820}). For a singular cardinal $\mu$, we introduced a positive consistency result under the assumption that $\mu$ is a limit of measurable cardinals (see \cite{MR2987137}).
An improvement of the method yielded the relation $\binom{\mu^+}{\mu}\rightarrow \binom{\mu^+}{\mu}^{1,1}_2$ for singulars which are limit of inaccessible cardinals, or even limit of strong limit cardinals. Of course, more properties are required for this consistency result, but one can force it. This appears in \cite{MR3201820}.

The process of replacing measurability by merely strong inaccessibility, excavated the fact that large cardinals are not really needed for these strong results. What we actually need is captured in the following metamathematical idea, which serves in many other theorems:

\begin{center}
\emph{Whenever GCH yields a negative result \\ an appropriate pcf assumption \\ may give a positive result}
\end{center}

The appropriate pcf assumption is not always at hand. But in our case it is possible to force the required assumption for every singular cardinal. Suppose $\mu$ is a singular cardinal. By forcing over a model of ZFC with the existence of a supercompact cardianl, we obtain a model of ZFC in which $\binom{\aleph_{\alpha+1}}{\aleph_\alpha} \rightarrow \binom{\aleph_{\alpha+1}}{\aleph_\alpha}^{1,1}_2$ holds for any prescribed limit non-inaccessible ordinal $\alpha$. This holds, in particular, for the first singular cardinal $\mu=\aleph_\omega$.

Before seeking into the detailed proofs, let us try to explain what we realy need for the strong polarized relation. Let $\mu$ be a singular cardinal, and let $\kappa$ be the cofinality of $\mu$. We may assume that $\mu$ is strong limit, and enough GCH holds below $\mu$. This can be forced for every $\mu$, and we shall need it in the proof. By a straightforward generalization of Theorem \ref{ErHaRado} we know that for getting the positive relation $\binom{\mu^+}{\mu}\rightarrow \binom{\mu^+}{\mu}^{1,1}_2$ we must increase $2^\mu$. In the case of countable cofinality we can do it with full GCH below $\mu$ (see \cite{MR0491183}, \cite{MR0491184} and \cite{MR823777}), and in the case of uncountable cofinality we will maintain enough GCH below $\mu$ to provide the combinatorial result on the one hand, and the ability to increase $2^\mu$ on the other hand.

Given a coloring $c:\mu^+\times\mu \rightarrow\theta$ we are trying to create a monochromatic product $A\times B$. We call $A$ the big component (its cardinality being $\mu^+$), and $B$ is the small component respectively. As the size of $B$ is the singular cardinal $\mu$, it is quite natural to compound $B$ from tiny pieces. So we fix a sequence of regular cardinals $\langle \mu_\varepsilon: \varepsilon<\kappa\rangle$ which tends to $\mu$. Our pcf assumptions will focus mainly on this sequence.

At first stage we deal with the small component $B$, so let us fix an ordinal $\alpha<\mu^+$, and concentrate in the restriction $c\upharpoonright \{\alpha\} \times\mu$. Moreover, we shall work with every $\mu_\varepsilon$ separately, so we concentrate in $c\upharpoonright \{\alpha\} \times\mu_\varepsilon$. For every $\gamma<\mu_\varepsilon$ we have a value $c(\alpha,\gamma)$, which is a color below $\theta$. Inasmuch as $\theta$ is small (and $\mu_\varepsilon$ is regular) we may pick a large set $B_\varepsilon$ of size $\mu_\varepsilon$ so that $c(\alpha,\gamma)$ is the same for every $\gamma\in B_\varepsilon$. Without loss of generality, the color is $0$.

The same process can be rendered for every $\varepsilon<\kappa$. Of course, the unified color for each $B_\varepsilon$ may vary from one $\varepsilon$ to another. But again, inasmuch as $\theta<\kappa=\cf(\kappa)$ we shall get $\kappa$-many $B_\varepsilon$-s with the same color. So we assume, without loss of generality, that the color is $0$ and it holds for each $\varepsilon<\kappa$. Set $B=\bigcup\{B_\varepsilon: \varepsilon<\kappa\}$, and the desired small component is at hand, with respect to our fixed $\alpha<\mu^+$.

We turn now to the big component, the set $A$. Fix for awhile some $\varepsilon<\kappa$. We can isolate some $B_\varepsilon\in[\mu_\varepsilon] ^{\mu_\varepsilon}$ for each $\alpha<\mu^+$, but notice that $B_\varepsilon$ depends on $\alpha$, so it would be better to call it $B^\alpha_\varepsilon$. Having different $B^\alpha_\varepsilon$-s (for different $\alpha$-s) is problematic, but $\mu^+$ is much larger than the number of possible $B^\alpha_\varepsilon$-s. Indeed, we have at most $\theta^{\mu_\varepsilon}$ options, and by the strong limitude of $\mu$ we know that $\theta^{\mu_\varepsilon}<\mu<\mu^+$. It follows that there is a set $A_\varepsilon$ of cardinality $\mu^+$ and a fixed set $B_\varepsilon$ of cardinality $\mu_\varepsilon$ so that $c\upharpoonright(A_\varepsilon \times B_\varepsilon)=0$.

There is only one step which is required to accomplish the construction. We have to collect all the $B_\varepsilon$-s. But in order to do so, we need to ensure that $A_\varepsilon$ is the same set for every $\varepsilon<\kappa$, or at least that the intersection of all the $A_\varepsilon$-s includes a set of size $\mu^+$. This is the last step, and it can never be accomplished in ZFC. In fact, if $2^\mu=\mu^+$ then the strong relation $\binom{\mu^+}{\mu}\rightarrow \binom{\mu^+}{\mu}^{1,1}_2$ fails. But let us try to understand what exactly prevents this last step.

In order to intersect $\kappa$-many sets, each of size $\mu^+$, and remain with a set of cardinality $\mu^+$, we need a $\kappa^+$-complete ultrafilter on $\mu^+$. Such an ultrafilter exists if some compact cardinal is situated between $\kappa$ and $\mu$. But in this direction we also need that each $A_\varepsilon$ is a set in the ultrafilter, and for this we need more and more completeness.

If we give up the full relation for the small component, then we can get the following result. Recall that $\kappa$ is compact iff the logic $\mathcal{L}_{\kappa\kappa}$ is compact. Compactness ensures the existence of $\kappa$-complete ultrafilters above $\kappa$.

\begin{claim}
\label{llimitofcompacts} Polarized relation for a limit of compact cardinals. \newline
Assume $\mu>\cf(\mu)=\kappa, \mu$ is a limit of compact cardinals. \newline
Then $\binom{\mu^+}{\mu}\rightarrow \binom{\mu^+}{\beta}^{1,1}_2$ for every $\beta<\mu$. \newline
Moreover, for every $\beta<\mu$ and every complete enough uniform ultrafilter $U_\beta$ on $\mu^+$ one can find a monochromatic product for which the big component belongs to $U_\beta$.
\end{claim}

\par\noindent\emph{Proof}. \newline
We shall use the notation of the discussion above this claim, so in particular a coloring $c$ is given, and all the cardinals and sets defined above.
Let $\beta$ be an ordinal below $\mu$. Let $\zeta<\kappa$ be the first ordinal such that $\beta<\mu_\zeta$, and denote $\theta^{\mu_\zeta}$ by $\tau$. Choose a $\tau^+$-complete ultrafilter $U_\beta$ on $\mu^+$. This is possible as $\mu$ is a limit of compact cardinals.

By the considerations above, one can choose $B_\varepsilon$-s and $A_\varepsilon$-s for every $\varepsilon<\zeta+1$, such that $A_\varepsilon\in U_\beta$ for every $\varepsilon<\zeta+1$. Set $A=\bigcap\{A_\varepsilon: \varepsilon<\zeta+1\}$, and $B=\bigcup\{B_\varepsilon: \varepsilon<\zeta+1\}$. It is easily verified that $A\times B$ is monochromatic for the coloring $c$, so the proof is complete.

\hfill \qedref{llimitofcompacts}

As clearly seen, this argument will never be true if we replace ``every $\beta<\mu$" by $\mu$ itself. Indeed, for $\mu$ we need a $\mu^+$-complete ultrafilter on $\mu^+$, which is impossible. One might suspect that this obstacle cannot be overcome, and the strong polarized relation (with respect to the pair $\mu,\mu^+$) must fail.

But here, in the last step of creating the monochromatic product, we can invoke pcf theory. Under some pcf assumptions we can show that the sets $A_\varepsilon$ are not arbitrary, and some properties are shared by every $\alpha<\mu^+$. This is exactly the missing part which is required to finish the construction of the full product.

The main idea is as follows. First we enumerate all the subsets of $\mu_\varepsilon$, for every $\varepsilon<\kappa$. It helps (although not essential) if $2^{\mu_\varepsilon}=\mu_\varepsilon^+$. Consider now a fixed ordinal $\alpha<\mu^+$. The sets of the form $B_\varepsilon$ appear in these enumerations, each $B_\varepsilon$ at some place of the pertinent enumeration. This draws a function $f_\alpha$ in the product $\prod\limits_{\varepsilon<\kappa} \mu_\varepsilon^+$. Repeating this process for every $\alpha<\mu^+$ we define $\mu^+$-many functions in the above product. Now pcf assumption enables us to find a single function $f$ in this product which bounds every $f_\alpha$. The assumption is that ${\rm tcf} (\prod\limits_{\varepsilon<\kappa} \mu_\varepsilon^+) \neq \mu^+$. Clearly, one has to increase $2^\mu$ above $\mu^+$ in order to get it, but nothing else is required.

How do we use this dominating function $f$? Well, it depends on the extra assumptions that we add. In the simplest case we assume that each $\mu_\varepsilon$ is a measurable cardinal (that was the idea behind the main theorem in \cite{MR2987137}). In this case we can arrange enough subsets of $\mu_\varepsilon$ in a $\subseteq^*$-decreasing chain, and just take the union of the $B_\varepsilon$-s which appear in the place $f(\varepsilon)$ for every $\varepsilon<\kappa$. In the more sophisticated theorems, culminating in the main theorem of the present paper, we iterate (finitely many times) such pcf arguments, and eliminate extra assumptions (like measurability) on the cardinals $\mu_\varepsilon$-s.

We emphasize that the corollaries from the main theorem below are, in a sense, best possible. The number of colors (which is $\theta$) must be below $\kappa$, the size of the monochromatic product cannot be larger of course, and the pcf assumptions can be forced for every singular $\mu$. We also indicate that the pair $(\mu, \mu^+)$ is just a special case, and the theorem below is phrased for the general case of $(\mu, \lambda)$ when $\lambda$ is any regular cardinal between $\mu$ and $2^\mu$. The strong polarized relation holds trivially for the pair $(\mu,\lambda)$ when $\lambda$ exceeds $2^\mu$, and the behaviour of a singular $\lambda$ is determined by that of $\cf(\lambda)$.

The paper consists of two parts. In the first part we list some open problems concerning the polarized relation.
We indicate that some progress has been made, see \cite{MR3129732} and also \cite{freesets} regarding the polarized relation and free sets.
In the second part (i.e., the last section) we provide a solution to some of them.

We thank the referee for many helpful suggestions.

\newpage

\section{Background}

Most of the problems in the list below are based on the notion of strong polarized relations. We divide them into two groups. One section is devoted to $\aleph_0$, and the second to uncountable cardinals. In this section we sketch knwon results, to be referred in the list of problems.

Polarized partition relations were introduced in \cite{MR0081864}, and investigated further in \cite{MR0202613}. Our starting point is the following negative relation of Erd\"os, Hajnal and Rado:

\begin{theorem}
\label{ErHaRado} If $2^{\aleph_0}=\aleph_1$ then $\binom{\omega_1}{\omega} \nrightarrow \binom{\omega_1}{\omega}^{1,1}_2$.
\end{theorem}

This negative relation is not a theorem of ZFC. First, it is consistent that $\binom{\omega_1}{\omega} \rightarrow \binom{\omega_1}{\omega}^{1,1}_2$ (e.g., if MA+$2^{\aleph_0}>\aleph_1$ holds, as proved in \cite{MR0371652}). Second, by enlarging $2^{\aleph_0}$ we can get a positive relation of the form $\binom{2^\omega}{\omega} \rightarrow \binom{2^\omega}{\omega}^{1,1}_2$. It was noted, first, in \cite{MR556894} in the specific case of $2^{\aleph_0}=\aleph_2$ (in the model of iterating Sacks reals) and proved in a more general setting in \cite{MR2927607} and \cite{MR3000439}:

\begin{theorem}
\label{ppositiveforcontinuum}
Suppose $\lambda>\aleph_1$ and $\lambda=\lambda^{\aleph_0}$ (hence $\cf(\lambda)>\aleph_0$).\newline
Then there is a forcing notion $\mathbb{P}$ so that $\mathfrak{c}=\lambda$ and $\binom{\lambda}{\omega} \rightarrow \binom{\lambda}{\omega}^{1,1}_2$ in ${\rm\bf V}^{\mathbb{P}}$.
\end{theorem}

A similar positive result holds upon replacing $\aleph_0$ by a supercompact cardinal $\mu$. This is also proved in \cite{MR2927607}:

\begin{theorem}
\label{ppositivecontinuumsupercompact} Suppose $\mu$ is a supercompact cardinal. \newline
Then it is consistent that $\binom{2^\mu}{\mu} \rightarrow \binom{2^\mu}{\mu}^{1,1}_2$.
\end{theorem}

A straightforward generalization of \ref{ErHaRado} asserts that $2^\kappa=\kappa^+$ implies $\binom{\kappa^+}{\kappa} \nrightarrow \binom{\kappa^+}{\kappa}^{1,1}_2$ for every infinite cardinal $\kappa$. As far as we know, a parallel positive result for $\kappa$ and $\kappa^+$ under the assumption $2^\kappa>\kappa^+$ is known only under assumptions of large cardinals. The following concerns regular limit cardinals (its proof is included in \cite{MR2927607}):

\begin{theorem}
\label{ppositivesupercompact} Suppose $\mu$ is a supercompact cardinal. \newline
Then in some forcing extension, $\mu$ remains supercompact and $\binom{\mu^+}{\mu} \rightarrow \binom{\mu^+}{\mu}^{1,1}_\theta$ holds true for every $\theta<\mu$.
\end{theorem}

A parallel result is proved in \cite{MR2987137} for singular cardinals:

\begin{theorem}
\label{ppositivesingular} Suppose there is a supercompact cardinal. \newline
Then in some forcing extension there is a singular cardinal $\mu$, a limit of measurable cardinals, so that $\binom{\mu^+}{\mu} \rightarrow \binom{\mu^+}{\mu}^{1,1}_\theta$ for every $\theta<\cf(\mu)$.
\end{theorem}

The next theorem concerns the cardinals below the splitting number. The proof of the countable case appears in \cite{MR2927607}, and the proof of the general case is situated in \cite{MR3201820}. For a thorough background on the splitting number $\mathfrak{s}$ see \cite{MR2768685}.
Here is the assertion:

\begin{theorem}
\label{bbelowssss} Suppose $\kappa<\lambda<\mathfrak{s}_\kappa$. \newline
The positive relation $\binom{\lambda}{\kappa} \rightarrow \binom{\lambda}{\kappa}^{1,1}_2$ holds iff $\cf(\lambda)\neq\kappa$.
\end{theorem}

The last theorem (from \cite{MR1606515}) to be mentioned in this section is almost strong. It gives a balanced form of the polarized relation, but the size of the big component can be $\alpha<\mu^+$ for every $\alpha<\mu^+$, and not $\mu^+$ itself:

\begin{theorem}
\label{sss586} Suppose $\mu$ is a limit of measurable cardinals. \newline
Then $\binom{\mu^+}{\mu} \rightarrow \binom{\alpha}{\mu}^{1,1}_2$ for every $\alpha<\mu^+$.
\end{theorem}

We state below some problems that remain unsolved in light of these results. Many of them are connected with cardinal invariants. We refer to \cite{MR2768685} for a comprehensive exposition of this subject. For a modern discussion on the polarized relation, see \cite{MR2768681}.

\newpage

\section{Countable problems}

Theorem \ref{ErHaRado} can be rephrased as follows: $\mathfrak{c}=\aleph_1$ implies $\binom{\mathfrak{c}}{\omega} \nrightarrow \binom{\mathfrak{c}}{\omega}^{1,1}_2$. We may ask whether $\mathfrak{c}$ is the correct cardinal invariant for this negative implication:

\begin{problem}
\label{c1} Suppose $\mathfrak{x}$ is a nicely defined invariant which satisfies $\mathfrak{x}=\aleph_1\Rightarrow \binom{\mathfrak{x}}{\omega} \nrightarrow \binom{\mathfrak{x}}{\omega}^{1,1}_2$. Does it follow that $\mathfrak{x}=\mathfrak{c}$?
\end{problem}

We are interested, in general, in the possibility that a cardinal invariant equals $\aleph_1$ yet a positive relation is consistent with respect to this invariant and $\aleph_0$. If $\mathfrak{x}=\aleph_1<\mathfrak{s}$ is consistent (e.g., $\mathfrak{x}=\mathfrak{a}$ or $\mathfrak{x}=\mathfrak{b}$) then $\mathfrak{x}=\aleph_1$ and $\binom{\mathfrak{x}}{\omega} \rightarrow \binom{\mathfrak{x}}{\omega}^{1,1}_2$ is consistent (see \cite{MR2927607}, Claim 2.4). So our problem applies to cardinal characteristics above the splitting number:

\begin{problem}
\label{c2} Is it consistent that $\mathfrak{s}=\aleph_1$ yet $\binom{\mathfrak{s}}{\omega} \rightarrow \binom{\mathfrak{s}}{\omega}^{1,1}_2$?
\end{problem}

The (consistency of the) positive relation $\binom{\mathfrak{s}}{\omega} \rightarrow \binom{\mathfrak{s}}{\omega}^{1,1}_2$ is established in \cite{MR3201820} in the case of $\mathfrak{s}=\aleph_2$. We may ask:

\begin{problem}
\label{c3} Is it consistent that $\mathfrak{s}>\aleph_2$ and $\binom{\mathfrak{s}}{\omega} \rightarrow \binom{\mathfrak{s}}{\omega}^{1,1}_2$?
\end{problem}

Another variant of this problem is related to the ultrafilter number. The positive result of \cite{MR3201820} comes from a model of Blass and Shelah (see \cite{MR879489}) in which $\mathfrak{u}<\mathfrak{s}=\mathfrak{c}$. Decreasing $\mathfrak{u}$ (or even $\mathfrak{r}$) is one way to get upward positive relations. But one may wonder if this is essential:

\begin{problem}
\label{c4} Is it consistent that $\mathfrak{s}=\mathfrak{u}$ and still $\binom{\mathfrak{s}}{\omega} \rightarrow \binom{\mathfrak{s}}{\omega}^{1,1}_2$?
\end{problem}

Or from the opposite point of view:

\begin{problem}
\label{c5} Is it consistent that $\mathfrak{s}<\mathfrak{c}$ and $\binom{\mathfrak{s}}{\omega} \rightarrow \binom{\mathfrak{s}}{\omega}^{1,1}_2$?
\end{problem}

Notice that a positive answer to Problem \ref{c2} entails a positive answer to the present problem.
As mentioned above, Problem \ref{c2} applies to every cardinal invariant above the splitting number. Our knowledge about the dominating number, in this respect, is slightly better. In particular, for every $\kappa=\cf(\kappa)>\aleph_1$ we can prove the consistency of $\mathfrak{d}= \kappa \wedge \binom{\mathfrak{d}}{\omega} \rightarrow \binom{\mathfrak{d}}{\omega}^{1,1}_2$. The proof is spelled out in \cite{MR3201820}, and the model comes from \cite{MR1005010}.
This gives rise to the following problems:

\begin{problem}
\label{c6} Is it consistent that $\mathfrak{d}=\aleph_1$ and $\binom{\mathfrak{d}}{\omega} \rightarrow \binom{\mathfrak{d}}{\omega}^{1,1}_2$?
\end{problem}

This problem is parallel to \ref{c2}. It is an open problem whether $\mathfrak{s}$ can be a singular cardinal. But it is known that $\mathfrak{d}$ might be singular. In fact, if $\aleph_0<\kappa=\cf(\kappa)\leq \cf(\lambda)\leq\lambda$ then $\mathfrak{b}=\kappa\wedge\mathfrak{d}=\lambda$ is consistent. This assertion is due to Hechler, see \cite{MR2768685} Theorem 2.5.
The model of Blass and Shelah from \cite{MR1005010} (in which $\mathfrak{u}<\mathfrak{d}$) requires the regularity of $\mathfrak{d}$. Hence the following may be phrased:

\begin{problem}
\label{c7} Is it consistent that $\mathfrak{d}$ is a singular cardinal and $\binom{\mathfrak{d}}{\omega} \rightarrow \binom{\mathfrak{d}}{\omega}^{1,1}_2$?
\end{problem}

More generally, we may try to eliminate the use of $\mathfrak{u}<\mathfrak{d}$ in order to get a positive relation. So the following (which is parallel to \ref{c4}) makes sense:

\begin{problem}
\label{c8} Is it consistent that $\mathfrak{d}\leq\mathfrak{u}$ and $\binom{\mathfrak{d}}{\omega} \rightarrow \binom{\mathfrak{d}}{\omega}^{1,1}_2$?
\end{problem}

And back again, we can try to look at this problem from the opposite direction of $\mathfrak{d}<\mathfrak{c}$. Here is the parallel to Problem \ref{c5}:

\begin{problem}
\label{c9} Assume $\binom{\mathfrak{d}}{\omega} \rightarrow \binom{\mathfrak{d}}{\omega}^{1,1}_2$, and one adds $\lambda$-many random reals (for some $\lambda>\mathfrak{d}$). Does the positive relation $\binom{\mathfrak{d}}{\omega} \rightarrow \binom{\mathfrak{d}}{\omega}^{1,1}_2$ still holds?
\end{problem}

The last two problems in this vein are about the independence number:

\begin{problem}
\label{c10} Is it consistent that $\mathfrak{i}=\aleph_1$ and $\binom{\mathfrak{i}}{\omega} \rightarrow \binom{\mathfrak{i}}{\omega}^{1,1}_2$?
\end{problem}

And with respect to the ultrafilter number:

\begin{problem}
\label{c11} Is it consistent that $\mathfrak{i}\leq\mathfrak{u}$ and $\binom{\mathfrak{i}}{\omega} \rightarrow \binom{\mathfrak{i}}{\omega}^{1,1}_2$?
\end{problem}

We indicate that the forcing construction from \cite{MR1175937}, in which the consistency of $\mathfrak{i}<\mathfrak{u}$ is established, might be helpful for some of the above problems.

The following problem is related to real-valued measurable cardinals. In the original construction of Solovay (see \cite{MR0290961}) it seems that such a cardinal carries a negative strong relation. But in \cite{MR2250550} appears a different way to introduce a real valued measurable cardinal in the universe. We ask:

\begin{problem}
\label{c12} Is the relation $\binom{\kappa}{\omega} \rightarrow \binom{\kappa}{\omega}^{1,1}_2$ consistent with $\kappa$ begin a real-valued measurable cardinal?
\end{problem}

Another direction of research concerns the spectrum of positive relations. Adding $\lambda$-many Cohen reals entails $\binom{\kappa}{\omega} \nrightarrow \binom{\kappa}{\omega}^{1,1}_2$ for every $\kappa\in[\aleph_1,\aleph]$ (an explicit proof appears in \cite{MR3000439}, Remark 2.4). The opposite situation is consistent as well. In \cite{MR3201820} it is proved that $\binom{\kappa}{\omega} \rightarrow \binom{\kappa}{\omega}^{1,1}_2$ for every $\kappa\in[\aleph_1,\aleph]$ is consistent in a model of $\mathfrak{c}=\aleph_2$. This invites the following:

\begin{problem}
\label{c13} Is it consistent that $\mathfrak{c}>\aleph_2$ and $\binom{\kappa}{\omega} \rightarrow \binom{\kappa}{\omega}^{1,1}_2$ for every $\kappa\in[\aleph_1,\aleph]$?
\end{problem}

This can be phrased in a slight stronger version. As mentioned above, if $\kappa<\mathfrak{s}$ and $\cf(\kappa)>\aleph_0$ then $\binom{\kappa}{\omega} \rightarrow \binom{\kappa}{\omega}^{1,1}_2$. It follows that if $\mathfrak{s}=\mathfrak{c}<\aleph_\omega$ then $\binom{\kappa}{\omega} \rightarrow \binom{\kappa}{\omega}^{1,1}_2$ for every $\kappa\in[\aleph_1,\aleph)$. We may ask:

\begin{problem}
\label{c14} Is it consistent that $\mathfrak{s}=\mathfrak{c}>\aleph_2$ and $\binom{\kappa}{\omega} \rightarrow \binom{\kappa}{\omega}^{1,1}_2$ for every $\kappa\in[\aleph_1,\aleph]$ (i.e., including $\kappa=\mathfrak{s}$ itself)?
\end{problem}

Of course, $\mathfrak{c}<\aleph_\omega$ in such models, as $\binom{\aleph_\omega}{\omega} \nrightarrow \binom{\aleph_\omega}{\omega}^{1,1}_2$. We say that $\aleph_\omega$ is a \emph{breaking point} in the positive spectrum of strong polarized relations. The following is a natural question:

\begin{problem}
\label{c15} Is it consistent that $\aleph_\omega$ is the only breaking point? More generally, can $\mathfrak{c}$ be arbitrarily large and the breaking points of the positive spectrum are just the cardinals of countable cofinality?
\end{problem}

By iterating Mathias forcing notions (relativised to some ultrafilter) we can get $\mathfrak{c}=\aleph_n$ and a positive spectrum which includes only one breaking point (see \cite{MR2927607}). One may ask if the positive spectrum can behave in a more complicated way:

\begin{problem}
\label{c16} Is it consistent that the polarized spectrum includes two breaking points below $\aleph_\omega$?
\end{problem}

Observe that if $\theta_0<\theta_1<\aleph_\omega$ are such points then $\mathfrak{s}\leq\theta_0$ and $\mathfrak{u}\geq\theta_1$. By and large, the splitting number is responsible (to a large extent) for the nature of this spectrum. We may wonder about invariants that might be located below $\mathfrak{s}$ and their relationship with $\mathfrak{s}$. The following is a typical question:

\begin{problem}
\label{c17} Is it consistent that $\mathfrak{a}=\mathfrak{s}$ and $\binom{\mathfrak{a}}{\omega} \rightarrow \binom{\mathfrak{a}}{\omega}^{1,1}_2$?
\end{problem}

The same can by phrased with respect to invariants that are always below the splitting number. It is proved in \cite{MR2444279} that $\binom{\mathfrak{p}}{\omega} \rightarrow \binom{\mathfrak{p}\ \alpha}{\omega\ \omega}^{1,1}_2$ for every $\alpha<\mathfrak{p}$. Due to Theorem \ref{bbelowssss}, if $\mathfrak{p}<\mathfrak{s}$ then $\binom{\mathfrak{p}}{\omega} \rightarrow \binom{\mathfrak{p}}{\omega}^{1,1}_2$. The following is natural:

\begin{problem}
\label{c18} Is it true that $\binom{\mathfrak{p}}{\omega} \rightarrow \binom{\mathfrak{p}}{\omega}^{1,1}_2$ iff $\mathfrak{p}<\mathfrak{s}$?
\end{problem}

A starting point might be the case of $\mathfrak{p}=\mathfrak{s}= \mathfrak{c}$. In the second stage, we may try to treat $\mathfrak{p}=\mathfrak{s}< \mathfrak{c}$ by collapsing $\mathfrak{c}$ using L\'evy collapse from \cite{MR0268037}. This invites a general problem:

\begin{problem}
\label{c19} Suppose $\aleph_1<\kappa<\lambda=\mathfrak{c}$ and $\binom{\kappa}{\omega} \rightarrow \binom{\kappa}{\omega}^{1,1}_2$. Let $\mathbb{P}={\rm Levy}(\kappa, \lambda)$ be the Levy collapse. Is it possible that $\binom{\kappa}{\omega} \nrightarrow \binom{\kappa}{\omega}^{1,1}_2$ in ${\rm \bf V}^{\mathbb{P}}$?
\end{problem}

Notice that if $\kappa=\aleph_1$ then the negative relation in ${\rm \bf V}^{\mathbb{P}}$ follows. But what about the slightly weaker relation:

\begin{problem}
\label{c20} Suppose $\kappa<\lambda=\mathfrak{c}$ and $\binom{\kappa}{\omega} \rightarrow \binom{\kappa}{\omega}^{1,1}_2$. Let $\mathbb{P}={\rm Levy}(\kappa, \lambda)$ be the Levy collapse. Is it possible that for some ordinal $\alpha<\kappa$ we have $\binom{\kappa}{\omega} \nrightarrow \binom{\kappa\ \alpha}{\omega\ \omega}^{1,1}_2$ in ${\rm \bf V}^{\mathbb{P}}$?
\end{problem}

We conclude this section with a problem of a different kind. It is proved in \cite{MR3000439} that the positive relation $\binom{\omega_1}{\omega} \rightarrow \binom{\omega_1}{\omega}^{1,1}_2$ is consistent with the weak diamond. We do not know what happens under the club principle:

\begin{problem}
\label{c21} Is it consistent that $\clubsuit_{\aleph_1}$ holds and $\binom{\omega_1}{\omega} \rightarrow \binom{\omega_1}{\omega}^{1,1}_2$?
\end{problem}

Observe that a positive answer requires $\clubsuit_{\aleph_1}$ with the negation of the continuum hypothesis. The consistency of this setting is proved in \cite{MR1623206}, but the forcing construction is a bit sophisticated. It suggests that a positive answer would not be easy.

\newpage

\section{Uncountable cardinals}

Recall that the consistency result $\binom{2^\omega}{\omega} \rightarrow \binom{2^\omega}{\omega}^{1,1}_2$ has been generalized in \cite{MR2927607} to the consistency of $\binom{2^\mu}{\mu} \rightarrow \binom{2^\mu}{\mu}^{1,1}_2$ for every supercompact cardinal (see Theorem \ref{ppositivecontinuumsupercompact} above). The proof is based on the generalized Mathias forcing. The generalization of Mathias forcing (relativised to some ultrafilter) requires the measurability of $\mu$. For iterating without destroying the measurability, we employ a Laver-indestructible supercompact $\mu$ (see \cite{MR0472529}). It feels that supercompactness is overwhelming, yet some large cardinals assumption is needed. So we phrase the following problems:

\begin{problem}
\label{d1} Does the negative relation $\binom{2^\kappa}{\kappa} \nrightarrow \binom{2^\kappa}{\kappa}^{1,1}_2$ hold for every successor cardinal $\kappa$?
\end{problem}

The same can be phrased for `small' large cardinals, in particular the ones below weakly compact:

\begin{problem}
\label{d2} Does the negative relation $\binom{2^\kappa}{\kappa} \nrightarrow \binom{2^\kappa}{\kappa}^{1,1}_2$ hold for every inaccessible but not weakly compact cardinal $\kappa$?
\end{problem}

In contrary, one may try to reduce the positive consistency result to a weaker assumption than supercompactness:

\begin{problem}
\label{d3} Let $\kappa$ be a measurable cardinal. Is there a forcing extension preserving the measurability of $\kappa$, in which
$\binom{2^\kappa}{\kappa} \rightarrow \binom{2^\kappa}{\kappa}^{1,1}_2$ hold? How about a strongly compact cardinal?
\end{problem}

Similar problems arise while replacing the pair $(\kappa,2^\kappa)$ by the pair $(\kappa,\kappa^+)$. Generally speaking, a positive strong relation of the form $\binom{\lambda}{\kappa} \rightarrow \binom{\lambda}{\kappa}^{1,1}_2$ becomes harder to achieve when $\lambda$ is closed to $\kappa$. For instance, if $\lambda=\kappa$ then a negative result follows for every infinite $\kappa$. Consequently, the possibility of $\binom{\kappa^+}{\kappa} \rightarrow \binom{\kappa^+}{\kappa}^{1,1}_2$ is an interesting problem.

By \cite{MR0371652}, Martin's Axiom implies $\binom{\aleph_1}{\aleph_0} \rightarrow \binom{\aleph_1}{\aleph_0}^{1,1}_2$. More generally, if $\mathfrak{s}>\aleph_1$ then $\binom{\aleph_1}{\aleph_0} \rightarrow \binom{\aleph_1}{\aleph_0}^{1,1}_2$ (as appears in Theorem \ref{bbelowssss}).
By \cite{MR3000439}, it is consistent that $\kappa$ is supercompact and $\mathfrak{s}_\kappa>\kappa^+$. Consequently, the relation $\binom{\kappa^+}{\kappa} \rightarrow \binom{\kappa^+}{\kappa}^{1,1}_2$ holds in this case.
The consistency of $\mathfrak{s}_\kappa>\kappa$ requires $\kappa$ to be a weakly compact cardinal (see \cite{MR1450512}), so the following is natural:

\begin{problem}
\label{d4} Is it consistent that $\binom{\kappa^+}{\kappa} \rightarrow \binom{\kappa^+}{\kappa}^{1,1}_2$ for a weakly compact cardinal $\kappa$ which is not measurable?
\end{problem}

And back again, to the opposite direction:

\begin{problem}
\label{d5} Suppose $\kappa=\cf(\kappa)$ is not a weakly compact cardinal. Is it provable in ZFC that $\binom{\kappa^+}{\kappa} \nrightarrow \binom{\kappa^+}{\kappa}^{1,1}_2$?
\end{problem}

The former problems deal with regular cardinals.
In the case of a singular cardinal $\lambda$ we can get a positive result according to Theorem \ref{ppositivesingular}. The proof requires that $\lambda$ is a limit of measurable cardinals. One may try to eliminate this assumption:

\begin{problem}
\label{d6} Let $\lambda$ be any strong limit singular cardinal. Is it consistent that $\binom{\lambda^+}{\lambda}\rightarrow \binom{\lambda^+}{\lambda}^{1,1}_2$? What about $\binom{2^\lambda}{\lambda}\rightarrow \binom{2^\lambda}{\lambda}^{1,1}_2$?
\end{problem}

Of course, a violation of the (local) continuum hypothesis on $\lambda$ must be employed. If $\lambda>\cf(\lambda)=\aleph_0$ this is possible, even if the GCH holds below $\lambda$. We can ask:

\begin{problem}
\label{d7} Suppose $\lambda>\cf(\lambda)=\aleph_0$, and $\theta<\lambda \Rightarrow 2^\theta=\theta^+$. Is it consistent that $\binom{\lambda^+}{\lambda}\rightarrow \binom{\lambda^+}{\lambda}^{1,1}_2$ or $\binom{2^\lambda}{\lambda}\rightarrow \binom{2^\lambda}{\lambda}^{1,1}_2$?
\end{problem}

Weaker results may hold without violating the GCH. Recall that if $\lambda$ is strong limit singular with $2^\lambda>\lambda^+$ then $\binom{\lambda^+}{\lambda}\rightarrow \binom{\lambda+1}{\lambda}^{1,1}_\theta$ for every $\theta<\cf(\lambda)$ (this is proved in \cite{MR1606515}). We can ask for the following strengthening:

\begin{problem}
\label{d8} Suppose $\lambda>\cf(\lambda)$, and $2^\lambda>\lambda^+$. Is it provable in ZFC that $\binom{\lambda^+}{\lambda}\rightarrow \binom{\alpha}{\lambda}^{1,1}_2$ for every $\alpha<\lambda^+$?
\end{problem}

Let us return to the case of a singular cardinal $\lambda$ being a limit of measurable cardinals. It is proved in \cite{MR1606515} that $\binom{\lambda^+}{\lambda}\rightarrow \binom{\alpha}{\lambda}^{1,1}_2$ for every $\alpha<\lambda^+$, in ZFC. As mentioned in the introduction, $\binom{\lambda^+}{\lambda}\rightarrow \binom{\lambda^+}{\lambda}^{1,1}_2$ fails to become a theorem of ZFC. But the intermediate situation remains open:

\begin{problem}
\label{d9} Suppose $\lambda>\cf(\lambda)$ is a limit of measurable cardinals. Is it provable in ZFC that $\binom{\lambda^+}{\lambda}\rightarrow \binom{\lambda^+\ \alpha}{\lambda\ \lambda}^{1,1}_2$ for every $\alpha<\lambda^+$?
\end{problem}

Even if this relation is not provable in ZFC, we may wonder if it is consistent with the local continuum hypothesis:

\begin{problem}
\label{d10} Suppose $\lambda>\cf(\lambda)$ is a limit of measurable cardinals. Is it consistent that $2^\lambda=\lambda^+$ and yet $\binom{\lambda^+}{\lambda}\rightarrow \binom{\lambda^+\ \alpha}{\lambda\ \lambda}^{1,1}_2$ for every $\alpha<\lambda^+$?
\end{problem}

On the other hand, one can get the impression that if $\lambda$ is not a strong limit cardinal then a positive consistency result is denied:

\begin{problem}
\label{d11} Is it provable in ZFC that $\binom{\lambda^+}{\lambda}\nrightarrow \binom{\lambda^+}{\lambda}^{1,1}_2$ whenever $\lambda>\cf(\lambda)$ and $\lambda$ is not strong limit?
\end{problem}

We conclude with a spectral problem:

\begin{problem}
\label{d12} Let $\lambda$ be a strong limit singular cardinal which is not a limit of measurable cardinals. Is it consistent that $\mu\in(\lambda,2^\lambda]\Rightarrow \binom{\mu}{\lambda}\rightarrow \binom{\mu}{\lambda}^{1,1}_2$?
\end{problem}

Under some extra assumptions on $\lambda$ we can give a positive answer (see \cite{MR3201820}). But the general case of a strong limit singular remains open.

\newpage

\section{Down to $\aleph_\omega$}

The theorem below is phrased in the general context of the pair $(\mu,\lambda)$. The specific cases of $\lambda=\mu^+$ and $\lambda=2^\mu$ are more interesting (see Problems \ref{d6} and \ref{d7}). We phrase the conclusions of this theorem, with respect to these cases.

\begin{theorem}
\label{mt} The main theorem. \newline
Assume that:
\begin{enumerate}
\item [$(a)$] $\mu>\cf(\mu)=\kappa$, and $\mu$ is a strong limit cardinal
\item [$(b)$] $\theta<\kappa$ ($\theta$ is the number of colors)
\item [$(c)$] $\mu<\lambda=\cf(\lambda)\leq 2^\mu$
\item [$(d)$] $\langle\mu_\varepsilon: \varepsilon<\kappa\rangle$ is an increasing sequence of regular cardinals
\item [$(e)$] $\mu=\bigcup\{\mu_\varepsilon: \varepsilon<\kappa\}$
\item [$(f)$] $\mu_\varepsilon=\sum\limits_{\alpha<\mu_\varepsilon} 2^{|\alpha|}$, for every $\varepsilon<\kappa$
\item [$(g)$] $2^{\mu_\varepsilon}=\mu_\varepsilon^+$ and $2^{\mu_\varepsilon^+}=\mu_\varepsilon^{++}$ for every $\varepsilon<\kappa$
\item [$(h)$] $\Upsilon_\ell={\rm tcf}(\prod\limits_{\varepsilon<\kappa} \mu_\varepsilon^{+\ell},J^{\rm bd}_\kappa)$ for $\ell\in\{0,1,2\}$
\item [$(i)$] $\lambda\notin\{\Upsilon_0,\Upsilon_1,\Upsilon_2\}$
\end{enumerate}
Then the strong polarized relation $\binom{\lambda}{\mu}\rightarrow \binom{\lambda}{\mu}^{1,1}_2$ holds.
\end{theorem}

\par\noindent\emph{Proof}.\newline
Let $c:\lambda\times\mu\rightarrow\theta$ be any coloring.
For every $\alpha<\lambda, \varepsilon<\kappa$ and $\iota<\theta$ we define the following set:

$$
A_{\alpha,\varepsilon,\iota}=\{\gamma<\mu_\varepsilon^+: c(\alpha,\gamma)=\iota\}
$$

We enumerate the members of $\mathcal{P}(\mu_\varepsilon^+)$ as $\langle B_{\varepsilon,i}: i<\mu_\varepsilon^{++}\rangle$. For every $\alpha<\lambda$ we define a function $f_{\alpha,\iota}\in \prod\limits_{\varepsilon<\kappa} \mu_\varepsilon^{++}$ by letting $f_{\alpha,\iota}(\varepsilon)$ be ${\rm min} \{i<\mu_\varepsilon^{++}: A_{\alpha,\varepsilon,\iota}= B_{\varepsilon,i}\}$. Now we define, for every $\alpha<\lambda$ the function $f_\alpha\in \prod\limits_{\varepsilon<\kappa} \mu_\varepsilon^{++}$ as follows:

$$
f_\alpha(\varepsilon)= {\rm sup}\{f_{\alpha,\iota}(\varepsilon): \iota<\theta\}
$$

Observe that $f_\alpha\in \prod\limits_{\varepsilon<\kappa} \mu_\varepsilon^{++}$ as $\theta$ is small and each $\mu_\varepsilon^{++}$ is regular. By assumption $(i)$, $\lambda\neq\Upsilon_2$, hence one can choose a function $f\in \prod\limits_{\varepsilon<\kappa} \mu_\varepsilon^{++}$ which bounds each $f_\alpha$. It means that $\alpha<\lambda\Rightarrow f_\alpha<_{J^{\rm bd}_\kappa}f$ for every $\alpha<\lambda$.

We may assume that $f(\varepsilon)\geq\mu_\varepsilon^+$ for every $\varepsilon<\kappa$.
By the definition of $J^{\rm bd}_\kappa$, for every $\alpha<\lambda$ there is $\varepsilon(\alpha)<\kappa$ for which $\varepsilon\geq\varepsilon(\alpha)\Rightarrow f_\alpha(\varepsilon)<f(\varepsilon)$. Since $\lambda=\cf(\lambda)>\kappa$, there exists a set $S\in[\lambda]^\lambda$ and an ordinal $\varepsilon<\kappa$ such that $\alpha\in S\Rightarrow \varepsilon(\alpha)=\varepsilon$.

All we need is a monochromatic product of size $\lambda\times\mu$, so the set $S$ and the union of sets of size $\mu_\delta$ for $\delta\in[\varepsilon,\kappa)$ suffice. Hence, we may assume (without loss of generality) that $S=\lambda$ and $\varepsilon=0$. This process of thinning out appears thrice throughout the proof, with the same mathematical meaning.

We focus on $\{B_{\varepsilon,i}: i<f(\varepsilon)\}$, which is a subcollection of $\mathcal{P}(\mu_\varepsilon^+)$, for every $\varepsilon<\kappa$. By re-enumerating its members as $\{B_{\varepsilon,i}^1: i<\mu_\varepsilon^+\}$ we can define (for every $\alpha<\lambda$ and $\iota<\theta$) a function $g_{\alpha,\iota}\in \prod\limits_{\varepsilon<\kappa} \mu_\varepsilon^+$ by $g_{\alpha,\iota}(\varepsilon)= {\rm min} \{i<\mu_\varepsilon^+: A_{\alpha,\varepsilon,\iota}= B_{\varepsilon,i}^1\}$. As before, for every $\alpha<\lambda$ we define:

$$
g_\alpha(\varepsilon)= {\rm sup}\{g_{\alpha,\iota}(\varepsilon): \iota<\theta\}
$$

Notice that $g_\alpha\in \prod\limits_{\varepsilon<\kappa} \mu_\varepsilon^+$ as well. We employ, again, assumption $(i)$ (this time we use $\lambda\neq \Upsilon_1$) to choose a bounding function $g\in \prod \limits_{\varepsilon<\kappa} \mu_\varepsilon^+$. In this round, we assume that $g(\varepsilon)\geq\mu_\varepsilon$ for every $\varepsilon<\kappa$, and without loss of generality, $g_\alpha(\varepsilon) <g(\varepsilon)$ for every $\alpha<\lambda$ and every $\varepsilon<\kappa$. This can be justified as in the above argument for the $f_\alpha$-s.

For the last round, we re-enumerate the collection $\{B_{\varepsilon,i}^1: i<g(\varepsilon)\}$ as $\{B_{\varepsilon,i}^2: i<\mu_\varepsilon\}$. We define $h_{\alpha,\iota}\in \prod\limits_{\varepsilon<\kappa} \mu_\varepsilon$ by $h_{\alpha,\iota}(\varepsilon)= {\rm min} \{i<\mu_\varepsilon: A_{\alpha,\varepsilon,\iota}= B_{\varepsilon,i}^2\}$. And back again, for every $\alpha<\lambda$ we define:

$$
h_\alpha(\varepsilon)= {\rm sup}\{h_{\alpha,\iota}(\varepsilon): \iota<\theta\}
$$

At this last stage we have $h_\alpha\in \prod\limits_{\varepsilon<\kappa} \mu_\varepsilon$. By the last part of $(i)$ we know that $\lambda\neq \Upsilon_0$, so a bounding function $h\in \prod \limits_{\varepsilon<\kappa} \mu_\varepsilon$ is at hand. Without loss of generality, it bounds each $h_\alpha$ for every $\varepsilon<\kappa$.

Now for each $\varepsilon<\kappa$ we define an equivalence relation $E_\varepsilon$ on the ordinals of $\mu_\varepsilon^+$ as follows:

$$
(\gamma_0 E_\varepsilon \gamma_1) \Leftrightarrow (\gamma_0\in B_{\varepsilon,j}^2 \Leftrightarrow \forall j<h(\varepsilon), \gamma_1\in B_{\varepsilon,j}^2)
$$

The number of equivalence classes is at most $2^{|h(\varepsilon)|} \leq\mu_\varepsilon< \mu_\varepsilon^+$. Choose (for every $\varepsilon<\kappa$) an equivalence class $X_\varepsilon$ of $E_\varepsilon$ of size $\mu_\varepsilon^+$. Observe that for every $\alpha<\lambda$ there is a color $\iota_{\alpha,\varepsilon}<\theta$ so that:

$$
\gamma\in X_\varepsilon \Rightarrow c(\alpha,\gamma)=\iota_{\alpha,\varepsilon}
$$

For every $\alpha<\lambda$ we choose a set $u_\alpha\in[\kappa]^\kappa$ and a fixed color $i_\alpha<\theta$ such that $\varepsilon\in u_\alpha \Rightarrow \iota_{\alpha,\varepsilon}=\iota_\alpha$. Since $2^\kappa<\lambda$ there is a set $A\in[\lambda]^\lambda$, a set $u\in[\kappa]^\kappa$ and a single color $\iota<\theta$ such that:

$$
\alpha\in A \wedge \varepsilon\in u \Rightarrow \iota_{\alpha,\varepsilon}=\iota
$$

Set $B=\bigcup\{X_\varepsilon: \varepsilon\in u\}$, so $|B|=\mu$. We claim that the product $A\times B$ is monochromatic for $c$. Indeed, if $\alpha\in A$ and $\beta\in B$ then there is an ordinal $\varepsilon\in u$ such that $\beta\in X_\varepsilon$, hence $c(\alpha,\beta)= \iota_{\alpha,\varepsilon}=\iota_\alpha=\iota$, and the proof is complete.

\hfill \qedref{mt}

The main theorem is phrased in a general fashion. We may derive some specific conclusions from it. These conclusions are based on known forcing constructions for singular cardinals, mainly from \cite{MR0491183},\cite{MR0491184} and \cite{MR823777}. We quote the following:

\begin{theorem}
\label{mmmagidor} The singular cardinals problem. \newline
Suppose $\mu$ is a supercompact cardinal, GCH holds above $\mu$ and $\delta<\mu_0<\mu$. \newline
There exists a forcing notion $\mathbb{P}$ which collapses cardinals only in the interval $(\mu_0,\mu]$, making $\mu$ a singular cardinal of cofinality $\kappa=\cf(\delta)$ and $\mu=\mu_0^{+\delta}$, such that in ${\rm\bf V}^{\mathbb{P}}$ the following hold:
\begin{enumerate}
\item [$(a)$] $\langle\mu_\varepsilon: \varepsilon<\kappa\rangle$ is an increasing sequence of regular cardinals
\item [$(b)$] $\mu=\bigcup\{\mu_\varepsilon: \varepsilon<\kappa\}$
\item [$(c)$] $\mu_\varepsilon=\sum\limits_{\alpha<\mu_\varepsilon} 2^{|\alpha|}$, for every $\varepsilon<\kappa$
\item [$(d)$] $2^{\mu_\varepsilon}=\mu_\varepsilon^+$ and $2^{\mu_\varepsilon^+}=\mu_\varepsilon^{++}$ for every $\varepsilon<\kappa$
\item [$(e)$] $\Upsilon_\ell={\rm tcf}(\prod\limits_{\varepsilon<\kappa} \mu_\varepsilon^{+\ell}, J^{\rm bd}_\kappa)$ for $\ell\in\{0,1,2\}$
\item [$(f)$] $\mu^+,2^\mu\notin\{\Upsilon_0,\Upsilon_1,\Upsilon_2\}$
\end{enumerate}
\end{theorem}

\par\noindent\emph{Proof}.\newline
By the constructions of \cite{MR0491183},\cite{MR0491184}, and mainly \cite{MR823777}. For a detailed proof we refer also to \cite{GaMgShF1211} and \cite{MR2768695}.

\hfill \qedref{mmmagidor}

Applying the above theorem to our problems, yields the following:

\begin{corollary}
\label{mmt} Strong polarized relations for singular cardinals. \newline
Suppose $\mu$ is a supercompact cardinal, GCH holds above $\mu$ and $\delta<\mu_0<\mu$. \newline
There is a forcing notion $\mathbb{P}$ as in the theorem above, making $\mu=\mu_0^{+\delta}$, such that:
\begin{enumerate}
\item [$(\aleph)$] The strong relation $\binom{\mu^+}{\mu}\rightarrow \binom{\mu^+}{\mu}^{1,1}_2$ holds in ${\rm\bf V}^{\mathbb{P}}$.
\item [$(\beth)$] The strong relation $\binom{2^\mu}{\mu}\rightarrow \binom{2^\mu}{\mu}^{1,1}_2$ holds in ${\rm\bf V}^{\mathbb{P}}$.
\end{enumerate}
\end{corollary}

\hfill \qedref{mmt}

Further corollaries, in this light, can be obtained (e.g., for getting an entirely positive spectrum of cardinals between $\mu^+$ and $2^\mu$, with respect to the strong polarized relation). We also indicate to an interesting consistency with GCH. By generalizing a remark of Foreman (quoted in \cite{MR1606515}), we can deduce that the consistency of the relation $\binom{\mu^+}{\mu}\rightarrow \binom{\mu^+}{\mu}^{1,1}_2$ under the negation of GCH implies the consistency of $\binom{\mu^+}{\mu}\rightarrow \binom{\alpha}{\mu}^{1,1}_2$ for every $\alpha<\mu^+$ with GCH. So this consistency holds, in light of the theorems above, for each singular cardinal.

\newpage

\bibliographystyle{amsplain}
\bibliography{arlist}

\end{document}